\documentclass[12pt]{article}

\usepackage{amsfonts}
\usepackage{xr}
\usepackage{amsmath,amsthm}
\usepackage{color,soul}
\usepackage[active]{srcltx}
\usepackage{hyperref}

\def\2{C^{1,2}(\R\times\R^N)}
\def\to{\rightarrow}
\def\e{\varepsilon}
\def\vp{\varphi}

\def\O{\Omega}
\def\Z{\mathbb{Z}}

\def\P{\mathbb{P}}
\def\R{\mathbb{R}}
\def\N{\mathbb{N}}
\def\tilde{\widetilde}
\def\.{\cdot}
\def\ul{\underline}
\def\ol{\overline}
\def\fa{\forall}
\def\t{\tilde}
\def\mc{\mathcal}
\def\supp{\text{\rm supp}\,}

\newcommand{\essinf}{\mathop{\rm ess{\,}inf}}

\newcommand{\be}{\begin{equation}}
\newcommand{\ee}{\end{equation}}
\newcommand{\baa}{\begin{array}}
\newcommand{\eaa}{\end{array}}
\newcommand{\ba}{\begin{eqnarray}}
\newcommand{\ea}{\end{eqnarray}}
\newcommand{\su}[2]{\genfrac{}{}{0pt}{}{#1}{#2}}

\newtheorem{thm}{\bf Theorem}[section]
\newtheorem{lem}[thm]{\bf Lemma}
\newtheorem{prop}[thm]{\bf Proposition}

\theoremstyle{definition}
\newtheorem{defi}[thm]{\bf Definition}
\newtheorem{rmq}[thm]{\bf Remark}
\newtheorem{ex}{\bf Example}
\newtheorem{hyp}[thm]{\bf Hypothesis}
\newtheorem{opb}{\bf Open problem}

\newenvironment{formula}[1]{\begin{equation}\label{eq:#1}}
                       {\end{equation}\noindent}
\def\Fi#1{\begin{formula}{#1}}
\def\Ff{\end{formula}\noindent}
\def\eq#1{{\rm(\ref{eq:#1})}}
\begin{document}
\date{}
\title{\bf{Propagation phenomena \\ for time heterogeneous \\ KPP
reaction-diffusion equations}}
\author{Gr\'egoire Nadin \thanks{CNRS, UMR 7598, Laboratoire Jacques-Louis
Lions, F-75005 Paris, France}  \and Luca Rossi \thanks{ Dipartimento di
Matematica P.~e A.,
  Universit\`a di Padova, via Trieste 63, 35121 Padova, Italy}}
\maketitle

\setstcolor{red}

\begin{abstract} We investigate in this paper propagation phenomena for the
heterogeneous reaction-diffusion equation  
$$ \partial_t u -\Delta u = f(t,u),\quad x\in\R^N,\ t\in\R,$$
where $f=f(t,u)$ is a KPP monostable nonlinearity which depends in a general way
on $t\in\R$. A typical $f$ which satisfies our hypotheses
is $f(t,u)=\mu (t) u(1-u)$, with 
$\mu\in L^\infty(\R)$ such that
$\essinf_{t\in\R}\mu (t)>0$. 
We first prove the existence of generalized transition waves (recently defined
in \cite{BerestyckiHamelgtw, Shenbistable})
for a given class of speeds.
As an application of this result, we obtain
the existence of random transition waves when $f$ is a random stationary
ergodic function with respect to $t\in\R$. Lastly, we prove some spreading
properties for the solution of the Cauchy problem.
\end{abstract}


\smallskip


\smallskip

\noindent This work was partially supported by the French ANR project {\em
Prefered}.\\
The second author was partially supported by GNAMPA - INdAM (Italy).

%%%%%%%%%%%%%%%%%%%%%%%%%%%%%%%%%%%%%%%%%%%%%%%%%%%%%%%%%

\section{Introduction}

This paper investigates the existence of traveling wave solutions for the heterogeneous reaction-diffusion equation 
\begin{equation}\label{eqprinc}
 \partial_t u -\Delta u = f(t,u),\quad x\in\R^N,\ t\in\R,
\end{equation}
where $f=f(t,u)$ vanishes when $u=0$ or $u=1$,
is strictly positive for $u\in(0,1)$ and is of KPP type, that is, 
$f(t,u)\leq f_u'(t,0)u$ for $(t,u)\in \R\times [0,1]$.

When $f$ is constant with respect to $t$, we recover the classical Fisher-KPP equation
\begin{equation}\label{eq-hom}
 \partial_t u -\Delta u = f(u), \quad x\in\R^N,\ t\in\R.
\end{equation}
It is well-known (see \cite{AronsonWeinberger, KPP}) that for all 
$c\geq c^*=2\sqrt{f'(0)}$, there exists a {\em planar traveling wave} 
of speed
$c$ in any direction  $e\in\mathbb{S}^{N-1}$, that is, a solution $u$ of
(\ref{eq-hom}) which can be written as 
$u(x,t) = \phi(x\cdot e -ct)$, with $\phi\geq 0$, $\phi(-\infty)=1$ and
$\phi(+\infty)=0$. In this case,  the profile $\phi$ of the planar traveling
wave satisfies the 
ordinary differential equation
$$-\phi''-c\phi' = f(\phi) \quad\hbox{in } \R.$$

When $f$ is periodic with respect to $t$, planar traveling waves no longer
exist and the relevant notion is that of 
{\em pulsating traveling wave}. Assume for instance 
that exists $T>0$ such that $f(t+T,u)= f(t,u)$ for all $(t,u) \in \R\times
[0,1]$. Then a pulsating traveling wave of speed 
$c$ in direction $e\in\mathbb{S}^{N-1}$ 
is a solution $u$ of (\ref{eqprinc}) which can be written as 
$u(x,t)=\phi (x\cdot e -ct,t)$, where $\phi\geq 0$, $\phi =\phi (z,t) $ is 
$T-$periodic with respect to
 $t$, $\phi (-\infty,t)=1$ and $\phi (+\infty,t)=0$ 
for all $t$. 
The profile $\phi$ then satisfies the time-periodic parabolic equation
$$\partial_t \phi-\partial_{zz}\phi-c\partial_z\phi = f(t,\phi)
\quad\text{in } \R\times\R.$$
The existence of pulsating traveling waves has been proved by Nolen, Rudd and
Xin \cite{NolenRuddXin}
 and the first author \cite{Nadinptf} under various hypotheses
in the more general framework of space-time periodic reaction-diffusion
equations
(see also \cite{Alikakos} 
for time periodic bistable equations). 
These results yield that, in the particular case of temporally periodic
equations, 
a pulsating traveling wave of speed $c$ exists if and only if $c\geq
c^*=2\sqrt{\langle \mu\rangle}$, where  
$\langle \mu\rangle=\frac{1}{T}\int_0^T f_u'(t,0)dt$. 

Note that the notion of pulsating traveling wave has first been introduced in
the framework of space periodic reaction-diffusion 
equations by Shigesada, Kawasaki and Teramoto
\cite{Shigesada} and Xin \cite{Xinperiodic} in parallel ways. Xin \cite{Xinperiodic}, Berestycki and Hamel \cite{BerestyckiHamel} 
and Berestycki, 
Hamel and Roques \cite{BerestyckiHamelRoques2} proved the existence of such
waves
in space periodic media under various hypotheses. In this case the minimal speed
$c^*$ for the existence of pulsating traveling waves is 
not determined using the mean average of $x\mapsto f_u'(x,0)$, but rather by
means of a family of periodic principal eigenvalues to characterize $c^*$. 

The case of a time almost periodic and bistable reaction term $f$ has been
investigated by Shen \cite{Shenap1, Shenap2}. A 
time heterogeneous nonlinearity is said to be {\em bistable} if 
there exists a smooth function $\theta : \R\to (0,1)$ such that 
\begin{equation}\label{f-bistable}
\left\{\begin{array}{l}
f(t,0)=f(t,1)=f(t,\theta (t))=0 \ \text{ for }t\in\R,\\
f(t,s)<0 \ \text{ for }t\in\R,\ s\in (0,\theta (t)),\\
f(t,s)>0 \ \text{ for }t\in\R,\  s\in (\theta (t),1).
       \end{array}
\right.
\end{equation}
Shen constructed examples where 
there exists no solution $u$ of the form $u(x,t)= \phi (x\cdot e -ct,t)$ such
that $\phi (-\infty,t)=1$, $\phi (+\infty,t)=0$ uniformly with respect to
$t\in\R$ and 
$\phi=\phi (z,t)$ is almost periodic with respect to $t$. She proved that the
appropriate notion of wave in almost periodic media incorporates a time
dependence
of the speed $c=c(t)$. Namely, she defined an {\em almost periodic traveling
wave} as a solution $u$ of (\ref{eqprinc}) which can be written as 
$u(x,t)= \phi (x\cdot e -\int_0^t c(s)ds,t)$, where $c=c(t)$ and $\phi=\phi
(z,t)$ are smooth functions which are almost periodic with respect to $t$ and 
$\phi (-\infty,t)=1$ and $\phi (+\infty,t)=0$ hold uniformly with respect
to $t\in\R$. 
She proved that for
nonlinearities $f$ which satisfy  (\ref{f-bistable}), there exists 
an almost periodic traveling
wave and that its profile and speed are uniquely determined up to translation of
the profile and to addition of the derivative of a time
almost periodic function to the speed. 

\smallskip

In order to handle general heterogeneous equations such as (\ref{eqprinc}), we will use 
the notion of (almost planar) generalized transition wave introduced by
Berestycki and Hamel \cite{BerestyckiHamelgtw} and 
Shen \cite{Shenbistable}. 
This definition is a natural extension of the earlier
notions. In particular, it enables 
a dependence of the speed with respect to time, as in almost periodic media.

\begin{defi}\label{def-gtf}
%1. 
An (almost planar) {\em generalized transition wave} in the direction
$e\in\mathbb{S}^{N-1}$ of equation (\ref{eqprinc}) is a time-global solution $u$
which can be written as  
$u(x,t)=\phi (x\cdot e -\int_0^t c(s)ds,t)$, where 
$c\in L^\infty(\R)$ is bounded and $\phi : \R\times \R \to [0,1]$
satisfies 
\begin{equation}
\lim_{z\to-\infty}\phi(z,t)=1 \ \hbox{ and }\ \lim_{z\to+\infty}\phi(z,t)=0\ \hbox{ uniformly in } t\in\R. 
\end{equation}
The functions $\phi$ and $c$ 
are respectively called the {\em profile} and the
{\em speed} of the generalized transition wave $u$.
\end{defi}

We will only consider almost planar waves in the present paper and thus we will
omit to mention the almost planarity in the sequel. Notice that the speed $c$ and the profile $\phi$ in Definition
\ref{def-gtf} are not univocally determined: they can be replaced by $c+\xi$
and $\phi (z+\int_0^t \xi(s)ds,t)$ respectively, for any $\xi\in L^\infty(\R)$
such
that $t\in \R\mapsto \int_0^t \xi (s)ds$ is bounded.
If $u$ is a generalized transition wave, then its profile $\phi$ satisfies
\begin{equation}\label{def-tf} 
\left\{\begin{array}{l}
\partial_t \phi-\partial_{zz}\phi -c(t)\partial_z \phi = f(t,\phi),\quad
z\in\R,\ t\in\R,\\ 
\displaystyle\lim_{z\to-\infty}\phi(z,t)=1 \ \hbox{ and }\
\lim_{z\to+\infty}\phi(z,t)=0\ \hbox{ uniformly in } t\in\R.\\ 
\end{array} \right.
\end{equation}

The existence of generalized transition waves has been proved by Shen in the
framework of time heterogeneous bistable  
reaction-diffusion equations \cite{Shenbistable}. Berestycki and Hamel \cite{BerestyckiHamelgtw2} also proved the 
existence of generalized transition waves 
for the monostable equation (\ref{eqprinc}) when $t\mapsto f(t,u)$ converges as
$t\to\pm\infty$ uniformly in $u$ 
(see Section \ref{sec:examples} below). 
Lastly, we learned while we were ending this paper that Shen proved the
existence of generalized transition waves when the coefficients are 
uniquely ergodic in a very recent paper \cite{Shenue}. We will describe more precisely the differences between Shen's approach and ours in 
Section \ref{sec:thm} below.
These are the only types of temporal heterogeneities for which the existence of generalized transition waves has been proved.

In space heterogeneous media, generalized transition waves have been shown to
exist for ignition type nonlinearity in dimension $1$. Namely, assuming that
$f=g(x)f_0(u)$,
with $g$ smooth, uniformly positive and bounded, and with $f_0$ satisfying
$$
\exists\theta\in (0,1),\qquad
\forall s\in [0,\theta]\cup \{1\},\quad f_0(s)=0,\qquad
\forall s\in (\theta,1),\quad f_0(s)>0.
$$
the existence of generalized transition waves has been obtained
in parallel ways by Nolen and Ryzhik \cite{NolenRyzhik} and Mellet, Roquejoffre
and Sire 
\cite{MelletRoquejoffreSire}. These generalized transition waves attract the solutions of the Cauchy problem associated with front-like initial data 
\cite{MelletNolenRoquejoffreRyzhik}. Zlatos gave some extensions of these results to multi-dimensional media 
and more general nonlinearities \cite{Zlatosdisordered}. 
Lastly, Nolen, Roquejoffre, Ryzhik and Zlatos constructed some space
heterogeneous KPP nonlinearities such that 
equation (\ref{eqprinc}) does not admit any generalized transition wave
\cite{NRRZ}.
A different generalization of the notion of traveling wave has been introduced
in the framework of bistable reaction-diffusion equations with a small 
perturbation in space of the homogeneous equation by Vakulenko and Volpert
\cite{VakulenkoVolpert}.

When the nonlinearity is homogeneous and of ignition or bistable type, we know
that there exists a unique speed associated with a planar traveling wave. The analogue
of this 
property in heterogeneous media is that, for bistable or ignition type nonlinearities, the generalized transition wave is unique up to translation in time 
(see \cite{BerestyckiHamelgtw2, MelletNolenRoquejoffreRyzhik, Shenbistable}). This means that the speed $t\mapsto c(t)$ of 
the generalized transition wave is unique in some sense. 

For monostable homogeneous nonlinearities, we know that there exists a unique
speed $c^*$ such that planar traveling waves of speed $c$ exist if and only if
$c\geq c^*$. 
Moreover, it is possible to construct solutions which behave as planar traveling
waves with different speeds when $t\to \pm\infty$ (see \cite{HamelNadirashvili})
and one can remark that these solutions are generalized transition waves.
Hence, 
we expect a wide range of speeds $t\mapsto c(t)$ to be associated with
generalized transition waves in heterogeneous media. 
The aims of the present paper are the following. 

\begin{itemize}
\item Prove the existence of generalized transition waves for time heterogeneous
monostable equations.

 \item Identify a set of speeds $t\mapsto c(t)$ associated with generalized transition waves.

\item Apply our results to particular nonlinearities such as random stationary ergodic. 
\end{itemize}

%%%%%%%%%%%%%%%%%%%%%%%%%%%%%%%%%%%%%%%%%%%%%%%%%%%%%%%%%%%%%%%%%%%%%%%%%%%%%%%%%%%%%%%%%%%%%%%%

\section{Statement of the results}

%%%%%%%%%%%%%%%%%%%%%%%%%%%%%%%%%%%%%%%%%%%%%%%%%%%%%%%%%%%%%%%%%%%%%%%%%%%%%%%%

\subsection{Hypotheses}\label{sec:hyp}

In this paper we just assume the nonlinear
term $f(t,u)$
to be bounded and measurable with respect to $t$. The notion of solution
considered is that of
strong solution: a subsolution (resp.~supersolution) $u$ of
(\ref{eqprinc}) is a function in $W^{2,1}_{N+1,loc}(\R^N\times\R)$ 
\footnote{ $W^{2,1}_{p}(\mc{Q})$, $\mc{Q}\subset\R^N\times\R$, stands for the
space of functions $u$ such that
$u,\,\partial_{x_{i}}u,\,\partial_{x_i x_j}u,\,\partial_t u\in L^p(\mc{Q})$.}
satisfying
$$\partial_t u -\Delta u\leq f(t,u)
\quad\text{(resp.~$\geq f(t,u)$)},\quad \text{for a.e.~}x\in\R^N,\ t\in\R,$$
and a solution is a function $u$ which is both a sub and a supersolution. 
It then follows from the standard parabolic theory that solutions belong to 
$W^{2,1}_{p,loc}(\R^N\times\R)$, for all $p<\infty$, and then
they are uniformly continuous by the embedding theorem.

The fundamental hypothesis we make on $f$ is that it is of KPP
type. 
Namely, $f(t,u)\leq\mu(t)u$, with $\mu(t):=f_u'(t,0)$, which means that
$f(t,\.)$ lies below its tangent at $0$ for all $t$. Hence,
we expect  
the linearization at $u=0$ to play a crucial role in
the dynamics of the equation. 

A typical $f$ we want to handle is 
\begin{equation} \label{examplef} f(t,u)=\mu (t) u(1-u), \end{equation} 
with
$\mu\in L^\infty (\R)$ and $\inf_{t\in\R}\mu(t)>0$. 
This function fulfills the following set of hypotheses, which are the ones
we will require in the statements of our results.

\begin{hyp}\label{hyp:KPP}
The function $f=f(t,u)$ 
satisfies $f(\.,u)\in L^\infty(\R)$, for all
$u\in[0,1]$, and, as a function
of $u$, is Lipschitz continuous in $[0,1]$ and of class $C^1$ in a neighborhood
of
$0$, uniformly with respect to $t\in\R$. Moreover,
setting $\mu(t):=
f_u'(t,0)$, the following properties hold:
\begin{equation}\label{hyp-pos} \text{for a.e.~} t\in\R,\quad
f(t,0)=f(t,1)=0,\qquad
\forall u\in(0,1),\quad \essinf_{t\in\R}f(t,u)>0,
\end{equation}
\begin{equation} \label{hyp-KPP}
\text{for a.e.~}(t,u)\in\R\times [0,1],\quad f(t,u)\leq \mu (t) u,
\end{equation}
\begin{equation} \label{hyp-lower}
\exists C>0,\ \gamma,\delta\in(0,1],\quad \text{for a.e.~}(t,u)\in\R\times
[0,\delta],\quad 
f(t,u)\geq \mu (t) u-Cu^{1+\gamma},
\end{equation}
\end{hyp}

Notice that (\ref{hyp-KPP}), (\ref{hyp-lower}) are fulfilled
if
$f(t,0)=0$ and $f(t,\.)$
is respectively concave and in
$C^{1+\gamma}([0,\delta])$, uniformly with respect to $t$.

%------------------------------------------------------------------

\subsection{Existence of generalized transition waves for general nonlinearities} \label{sec:thm}

Since the range of speeds associated with planar traveling
waves in homogeneous media is given by $[2\sqrt{f'(0)},+\infty)$, we expect
to find similar constraints on the 
speeds of generalized transition waves in heterogeneous media. The
constraint we will exhibit depends on the {\em least mean} of the speed. 

\begin{defi}\label{def-mean}
For any given function $g\in L^\infty(\R)$, we define
$$\underline{g}:=\sup_{T>0}\,\inf_{t\in\R}
\frac{1}{T}\int_t^{t+T}g(s)\,ds$$
and we call this quantity the {\em least mean} of the function $g$
(over $\R$). 
\end{defi}

The definition of least mean does not change if one
replaces the $\sup_{T>0}$ with $\lim_{T\to+\infty}$ (see Proposition
\ref{reformulation} below). Hence, if $g$ admits a {\em mean value}
$\langle
g\rangle$, i.e., if there 
exists 
\Fi{mv}
\langle g\rangle:=\lim_{T\to+\infty}
\frac{1}{T}\int_t^{t+T}g(s)\,ds,\quad\text{uniformly with respect to }t\in\R,
\Ff
then $\ul g=\langle g\rangle$.

Notice that (\ref{hyp-pos}) and (\ref{hyp-KPP}) yield $\ul \mu>0$.
We are now in position to state our main result. 

\begin{thm}\label{thm-existence}
Assume that $f$ satisfies Hypothesis \ref{hyp:KPP}
and let $e\in\mathbb{S}^{N-1}$.

1) For all $\gamma>2\sqrt{\ul \mu}$,  
there exists a generalized transition wave $u$ in direction $e$ with a speed $c$
such that $\ul c= \gamma$ 
and a profile $\phi$ which is decreasing with respect to $z$.

2) There exists no generalized transition wave $u$ in direction $e$
with a speed $c$ such that 
$\ul c< 2\sqrt{\ul\mu}$.
\end{thm}

The speeds we construct have the particular form
\begin{equation} \label{formc}
 c(t)=\frac{\mu(t)}{\kappa}+\kappa, \qquad \kappa\in (0, \sqrt{\ul \mu}).
\end{equation}
Hence, they keep some properties of the
function $\mu$. 
In particular, if $\mu$ admits a mean value $\langle \mu \rangle$ then
our result implies that a generalized transition wave with a
speed $c$ such that $\langle c \rangle=\gamma$ exists if $\gamma>2\sqrt{\langle \mu
\rangle}$ and does not exist if $\gamma<2\sqrt{\langle \mu \rangle}$.
Of course, this construction is not exhaustive:
there might exist generalized transition waves with a speed $c$ which cannot be
written in the form (\ref{formc}), as exhibited in Example \ref{E:2speeds}
below. 
More generally, trying to characterize the set of 
speeds associated with generalized transition waves is a very hard task, since
this notion covers many types of speeds
(see \cite{BerestyckiHamelgtw, HamelNadirashvili}). Indeed, it is still an open
problem even in the case of homogeneous equations.

We did not manage to prove the existence of generalized transition waves with a
speed with least mean $2\sqrt{\ul \mu}$. We leave this extension as an open
problem, that we will discuss in Examples \ref{E:speeds} and \ref{E:minimal} of
Section \ref{sec:examples} below. 

Theorem \ref{thm-existence} shows that the range of least means of the
speeds associated with generalized transition waves is a half-line, with
infimum $2\sqrt{\ul\mu}$.
If instead one considers other notions of mean then
the picture is far from being complete: our existence result implies that, for
every notion of mean $\mc{M}$ 
satisfying 
$$\fa g\in L^\infty(\R),\quad
\mc{M}(g)\geq\ul g,\qquad
\fa\alpha,\beta>0,\quad\mc{M}(\alpha g+\beta)=\alpha\mc{M}(g)+\beta,$$
and every $\kappa\in(0,\sqrt{\ul
\mu})$, there exists a wave with speed $c$ satisfying
$\mc{M}(c)=\frac{\mc{M}(\mu)}\kappa+\kappa$, whereas there are no waves with
$\mc{M}(c)<2\sqrt{\ul\mu}$. But if $\mc{M}(\mu)>\ul\mu$ then
$2\sqrt{\ul\mu}<\frac{\mc{M}(\mu)}{\sqrt{\ul\mu}}+\sqrt{\ul\mu}$,
whence there is a gap between the thresholds of the existence and non-existence
results.

\bigskip

In order to conclude this section, we briefly comment the differences between
Shen's approach in \cite{Shenue} and the present one. First, Shen only
considered 
uniquely ergodic coefficients. We refer to \cite{Shenue} for a precise definition but we would like to point out that if $\mu$ is uniquely ergodic, 
then $\langle \mu \rangle=\lim_{T\to +\infty} \frac{1}{T}\int_{t}^{t+T}\mu (s)ds$ exists uniformly with respect to $t\in\R$.  
This hypothesis is quite restrictive since it excludes, for example, general random stationary ergodic coefficients. Under this hypothesis, Shen proves 
that a
generalized transition wave with a uniquely ergodic speed $c=c(t)$ satisfying
$\langle c \rangle = \gamma$ exists if and only if
$\gamma\geq 2\sqrt{\langle \mu\rangle}$. This is a slightly stronger result than
ours since it provides existence in the critical case $\gamma= 2\sqrt{\langle
\mu \rangle}$. Lastly, 
Shen's approach is different since she uses a dynamical systems setting, while
we use a PDE approach inspired by \cite{BerestyckiHamel}. 

%----------------------------------------------------------------------------

\subsection{Application to random stationary ergodic equations}

We consider the reaction-diffusion equation with random nonlinear term
\begin{equation}\label{eq-random}
 \partial_t u -\Delta u = f(t,\omega,u),\quad
x\in\R^N,\ t\in\R,\ \omega\in\Omega.
\end{equation} 
The function $f:\R\times \O\times [0,1]\to \R$ is a random function defined on
a probability space $(\O,\P,\mathcal{F})$. We assume that $(t,u)\mapsto
f(t,\omega,u)$
satisfies Hypothesis \ref{hyp:KPP} for almost every
$\omega\in\O$, and in addition that
\begin{equation} \label{hyp-vois1}
\fa t\in\R,\quad u\mapsto f(t,\omega,u)/u \hbox{ is nonincreasing in }
[1-\delta,1],
\end{equation}
where the constants $\gamma$ and $\delta$ depend on $\omega$.
Notice that (\ref{hyp-vois1}) 
is satisfied in particular when $u\mapsto f(t,\omega,u)$
is nonincreasing in $(1-\delta,1)$, since $f$ is positive there.
We further suppose that $f(t,\omega,u)$ is a stationary ergodic random
function with respect to $t$. Namely, there exists a group $(\pi_t)_{t\in\R}$ of
measure-preserving transformations of $\O$
such that 
$$\forall (t,s,\omega,u)\in\R\times
\R\times \O\times [0,1],\quad
f(t+s,\omega,u)=f(t,\pi_s \omega,u),$$ 
and for all $A\in\mathcal{F}$, if $\pi_t A=A$ for all $t\in\mathbb{R}$, then
$\P(A)=0$ or $1$.

A generalization of the notion of traveling waves for equation (\ref{eq-random})
has been given by Shen in \cite{Shenrandom}.

\begin{defi}\label{def-randomtw} {\bf (see \cite{Shenrandom}, Def. 2.3)}
A {\em random transition wave} in the direction $e\in\mathbb{S}^{N-1}$ of
equation (\ref{eq-random}) is a
function $u:\R^N\times \R \times \Omega \to [0,1]$
which satisfies:
\begin{itemize}
 \item There exist two 
bounded
measurable functions $\tilde{c}:\Omega \to \R$ and
$\tilde{\phi}:\R \times\Omega\to[0,1]$ such that $u$ can be written as
$$u(x,t,\omega) = \tilde{\phi}(x\cdot e -\int_0^t \tilde{c}(\pi_s
\omega)ds,\pi_t \omega) \hbox{ for all } (x,t,\omega)\in \R^N\times
\R\times\Omega.$$
\item For almost every $\omega\in\Omega$, $(x,t)\mapsto u(x,t,\omega)$ is
a solution of (\ref{eq-random}).
\item For almost every
$\omega\in\Omega$, $\displaystyle\lim_{z\to-\infty}\tilde{\phi}(z,\omega)=1$ and
$\displaystyle\lim_{z\to+\infty}\tilde{\phi}(z,\omega)=0$.
\end{itemize}
The functions $\tilde{\phi}$ and $\tilde{c}$ are respectively called the
{\em random profile} and the {\em random speed} of the random transition
wave $u$.
\end{defi}

Notice that if (\ref{eq-random}) admits a generalized transition wave for
a.e.~$\omega\in\O$, and the associated profiles $\phi(z,t,\omega)$ and speeds
$c(t,\omega)$ are stationary ergodic with respect to $t$, then the functions
$\t\phi(z,\omega):=\phi(z,0,\omega)$ and $\t c(\omega):=c(0,\omega)$ are the
profile and the speed of a random transition wave.

The existence of random transition waves has been proved in the framework of
space-time random stationary ergodic bistable nonlinearities by Shen
\cite{Shenrandom} and in the 
framework of space random stationary ergodic ignition type nonlinearities by
Nolen and Ryzhik \cite{NolenRyzhik} (see also \cite{Zlatosdisordered} for some
extensions). 

Starting from Theorem \ref{thm-existence}, we are able to characterize the
existence of random transition waves in terms of the least mean of their speed.
For a stationary ergodic function $g:\R\times\O\to\R$, the least mean of
$t\mapsto g(t,\omega)$ is independent of $\omega$, for every $\omega$ in a set
of probability $1$ (see Proposition \ref{pro:lmergodic} below). We call this
quantity the least mean of the random function $g$, and we denote it by $\ul g$.

\begin{thm}\label{thm-existencerandom}
Let $e\in\mathbb{S}^{N-1}$. Under the previous hypotheses,
for all $\gamma>2\sqrt{\ul\mu}$,  
there exists a random transition wave $u$ in direction $e$ with random speed $\t
c$ such that $c(t,\omega):=\tilde{c}(\pi_t\omega)$ has least mean $\gamma$,
and a random profile $\t\phi$ which is decreasing with respect to $z$.
\end{thm}

This result is not an immediate corollary of Theorem \ref{thm-existence}. 
In
fact, for given $\kappa\in (0,\sqrt{\ul \mu})$ and almost every $\omega\in\O$,
Theorem \ref{thm-existence} provides a generalized transition wave  
with speed $c(\cdot,\omega)$ satisfying 
$c=\frac{\mu}{\kappa}+\kappa$. 
Hence, $c$ is stationary ergodic in $t$, but it is far from being obvious that
the same is true for the profile $\phi$.
Actually, to prove this we require the additional hypothesis
(\ref{hyp-vois1}). Instead, the non-existence of random transition waves with
speeds $c$ satisfying $\ul c<2\sqrt{\ul\mu}$ follows directly from Theorem
\ref{thm-existence}.

%%%%%%%%%%%%%%%%%%%%%%%%%%%%%%%%%%%%%%%%%%%%%%%%%%%%%%%%%%%%%%%%%%%%%%%%%%%%%%%%

\subsection{Spreading properties}\label{sec:spreading}

When the nonlinearity $f$ does not depend on $t$, Aronson and Weinberger
\cite{AronsonWeinberger} proved that  
if $u$ is the solution of the associated Cauchy problem, with an initial datum
which is ``front-like'' in direction $e$, then  
for all $\sigma>0$, 
$$\lim_{t\to +\infty}\,\inf_{x\leq
(2\sqrt{f'(0)}-\sigma)t}u(x,t)=1,\qquad
\lim_{t\to +\infty}\, \inf_{x\geq
(2\sqrt{f'(0)}+\sigma)t}u(x,t)=0.$$
This result is called a {\em spreading property} and means that the level-lines
of $u(t,\cdot)$ behave like $2\sqrt{f'(0)}t$ as $t\to +\infty$. 
The aim of this section is to extend this property to the Cauchy problem
associated with (\ref{eqprinc}), namely,
\Fi{Cauchy}
\left\{\begin{array}{ll}
      \partial_t u -\Delta u = f(t,u),& x\in\R^N,\ t>0,\\
      u(x,0)=u_0(x), & x\in\R^N.
      \end{array}\right.
\Ff
Our result will involve once again the least mean of $\mu$, but this time over $\R_+$, because the
equation is defined only for $t>0$. For a given function $g\in L^\infty
((0,+\infty))$, we set
$$
\ul g_+:=\sup_{T>0}\,\inf_{t>0}
\frac{1}{T}\int_t^{t+T}g(s)\,ds.
$$
We similarly define the upper mean $\ol g^+$:
$$
\ol g^+:=\inf_{T>0}\,\sup_{t>0}
\frac{1}{T}\int_t^{t+T}g(s)\,ds.
$$
In \cite{BHN} Berestycki, Hamel and the first author partially extended the
result of \cite{AronsonWeinberger} to
general space-time heterogeneous equations. They showed in particular
that the level-lines  
of $u(t,\cdot)$ do not grow linearly and can oscillate. They obtained some
estimates on the location of these level-lines, which are
optimal
when $t\mapsto f(t,u)$ 
converges as $t\to +\infty$ locally in $u$, but not when $f$ is periodic for
example. 
These properties have been improved by Berestycki and the first author in
\cite{BN}, by using the notion of {\em generalized principal eigenvalues} in
order to estimate  
more precisely the maximal and the minimal linear growths of the location of the
level-lines of $u(t,\cdot)$.
When $f$ only depends on $t$, as in the present paper, they proved that if 
$u_0\in C^0(\R^N)$ is such that $0\leq u_0\leq 1$, $u_0\not\equiv 0$ and it is
compactly supported, then the solution $u$ of \eq{Cauchy} satisfies
\Fi{01}
\fa e\in\mathbb{S}^{N-1},\quad\lim_{t\to+\infty}u(x+\gamma te,t)=
\left\{\begin{array}{lll}
1 &\hbox{ if }& 0\leq\gamma <2\sqrt{\ul\mu_+}\\
0 &\hbox{ if }& \gamma>2\sqrt{\ol\mu^+}\\
\end{array}\right. 
\quad\hbox{locally in } x.
\Ff
In \cite{BN}, this result follows from a more general statement, proved in the
framework of space-time heterogeneous equations using homogenization
techniques. 
Here we improve \eq{01} by decreasing the threshold for the convergence to 0.
Our proof is based on direct arguments.

\begin{prop}\label{prop-spreading} Assume that $f$ satisfies Hypothesis
\ref{hyp:KPP} and let $u_0\in C^0 (\R^N)$ be such that $0\leq u_0 \leq
1$, $u_0\not\equiv0$. Then the solution $u$ of \eq{Cauchy} satisfies
\begin{equation}\label{eq:spreadinginf}
\fa\gamma<2\sqrt{\ul\mu_+},\quad
\lim_{t\to +\infty}\inf_{|x|\leq \gamma t} u(x,t)=1.
\end{equation}
If in addition $u_0$ is compactly supported then 
\begin{equation}\label{eq:spreadingsup}
\fa\sigma>0,\quad
\lim_{t\to+\infty}\sup_{|x|\geq 2\sqrt{t\int_0^t \mu (s)ds} +\sigma t} u(x,t)
=0.
\end{equation}
\end{prop}

If $u_0$ is ``front-like'' in the direction $e$, then 
$|x|\geq 2\sqrt{t\int_0^t \mu (s)ds} +\sigma t$
can be replaced by
$x\cdot e\geq 2\sqrt{t\int_0^t \mu (s)ds} +\sigma t$ in \eq{spreadingsup}.

\begin{rmq} Proposition \ref{prop-spreading} still
holds if (\ref{hyp-lower}),
(\ref{hyp-vois1}) are not satisfied and, in case of \eq{spreadinginf},
the KPP condition (\ref{hyp-KPP}) can also be dropped. 
\end{rmq}

If $\frac1t\int_0^t \mu (s)ds\to\ul\mu_+$ as $t$ goes to $+\infty$ then the result of Proposition \ref{prop-spreading} is optimal. 
Otherwise it does not describe in an exhaustive way the large time behavior of $u$.

\begin{opb} Assume that the hypotheses of Proposition \ref{prop-spreading}
hold and that
$$\ul\mu_+<\gamma<\liminf_{t\to+\infty}\frac1t\int_0^t \mu (s)ds.$$
What can we say about $\lim_{t\to+\infty}u(\gamma t e,t)$.
\end{opb}

%--------------------------------------------------------------------------

\subsection{Examples.} \label{sec:examples}

We now present some examples in order to illustrate the notion of generalized
transition waves and to discuss the optimality of our results. 

\begin{ex}\label{E:means}
{\em Functions without uniform mean.}\\
Set $t_1:=2$ and, for $n\in\N$,
$$\sigma_n:=t_n+n,\qquad \tau_n:=\sigma_n+n,\qquad t_{n+1}:=\tau_n+2^n.$$
The function $\mu$ defined by
$$\mu(t):=\left\{\begin{array}{ll}
                  3 & \text{if }t_n<t<\sigma_n,\ n\in\N,\\
		  1 & \text{if }\sigma_n<t<\tau_n,\ n\in\N,\\
                  2 & \text{otherwise}
                 \end{array}
\right.$$
satisfies
$$\ul\mu_+=1<\lim_{t\to+\infty}\frac1t\int_0^t \mu (s)ds=2<\ol\mu_+=3.$$
Therefore, $\mu$ does not admit a uniform mean $\langle\mu\rangle$
(over $\R_+$).
\end{ex}

\begin{ex}\label{E:speeds}
{\em Generalized transition waves with various choices of speeds.}\\
Fix $\alpha>1$ and consider the homogeneous reaction-diffusion equation
$$
\partial_t u -\partial_{xx}u = u(1-u)(u+\alpha).
$$
A straightforward computation yields that the function
$$U (z) = \frac{1}{1+e^{z/\sqrt2}}$$
is the profile of a planar traveling wave of speed $c_0 =
\sqrt{2} \alpha +\frac{1}{\sqrt{2}}$. Note that $c_0$ is strictly larger than the minimal speed $2\sqrt{\alpha}$ 
for the existence of traveling waves.

We now perturb this equation by adding some time heterogeneous bounded
function $\xi\in C^0(\R)$ in the nonlinearity:
\begin{equation}\label{example-2}
\partial_t v -\partial_{xx}v= v(1-v)(v+\alpha-\frac{1}{\sqrt{2}} \xi (t)),
\end{equation}
with $\|\xi\|_\infty < \sqrt{2}(\alpha-1)$ so that this equation is still of KPP
type. 
Let $$v(x,t):=U(x-c_0 t +\int_0^t \xi (s)ds).$$ 
Since $U' = -\frac{U}{\sqrt{2}}
(1-U)$, one readily checks that $v$ is a generalized transition wave of
equation (\ref{example-2}), with
speed 
$c(t)= c_0 +\xi (t)$. If $\langle\xi \rangle =0$,
then this generalized transition wave has a global mean speed $c_0$. But as
$\xi$ 
is arbitrary, the fluctuations around this global mean speed can be large. 
This example shows that the speeds associated with generalized transition waves
can be very general, depending on the structure of the heterogeneous equation. 
\end{ex}

\begin{ex}\label{E:minimal}
{\em Generalized transition waves with speeds $c$ satisfying $\ul c =
2\sqrt{\ul \mu}$.}
We can generalize the method used in the previous example to obtain generalized
transition waves with a speed with minimal least mean. Consider any homogeneous
function
$f: [0,1]\to [0,\infty)$ such that $f$ is Lipschitz-continuous, $f(0)=f(1)=0$,
$f(s)>0$ if $s\in (0,1)$ and $s\mapsto f(s)/s$ is nonincreasing. 
Then we know from \cite{AronsonWeinberger} that for all $c\geq 2\sqrt{f'(0)}$,
there exists a decreasing function $U_c\in C^2(\R)$ such that $U_c(-\infty)=1$, 
$U_c (+\infty)=0$ and $-U_c''-cU_c'=f(U_c)$ in $\R$. 
It is well-known that $-U_c''(x)/U_c'(x)\to \lambda_c$ as $x\to
+\infty$, where $\lambda_c>0$ is the smallest root of
$\lambda\mapsto-\lambda^2+c\lambda_c-f'(0)$. 
Moreover, writing the equation satisfied by $U_c'/U_c$, one can prove that
$U_c'\geq -\lambda_c U_c$ in $\R$.
It follows that the function $P_c
:  [0,1] \to[0,+\infty)$ defined by $P_c(u):=-U_c'(U_c^{-1}(u))$ for $u\in(0,1)$
and $P_c (0)=P_c (1)=0$ is Lipschitz-continuous.
Furthermore, it is of KPP type, because
$$P_c'(0)=-\lim_{x\to+\infty}\frac{U_c''(x)}{U_c'(x)}=\lambda_c,
\qquad \frac{P_c(u)}u=-\frac{U_c'(U_c^{-1}(u))}{U_c(U_c^{-1}(u))}
\geq\lambda_c.$$
We now consider a given function $\xi \in C^\alpha(\R)$ with least mean $\ul
\xi=0$. The function $v$ defined by $v(x,t):= U_c (x-ct -\int_0^t \xi (s)ds)$
satisfies
$$\partial_t v -\partial_{xx} v = f(v) +\xi (t) P_c(v)=: g(t,v)\quad \hbox{in }
\R.$$
It is clearly a generalized transition wave of this equation, with speed $c_1
(t)= c+\xi (t)$. 
Let now see what Theorem \ref{thm-existence} gives. Here, $\mu (t)=
g_u'(t,0)=f'(0)+\lambda_c \xi (t)$ and thus $\ul \mu = f'(0)$ since $\ul
\xi=0$. 

If $c>2\sqrt{f'(0)}$, then we know that there exists 
a generalized transition wave $w$ with a speed $c_2(t)=\frac{\mu
(t)}{\kappa}+\kappa$, for some $0<\kappa<\sqrt{f'(0)}$ such that $\ul c_2=
c$. These two
conditions impose 
$\kappa=\lambda_c$ and thus $c_2 (t)=\frac{f'(0)}{\lambda_c}+\xi (t)+\lambda_c=
c_1 (t)$. Hence, $c_1\equiv c_2$, which means that the speed obtained through
Theorem
\ref{thm-existence} is the speed of the generalized transition wave $v$. 

The case $c=2\sqrt{f'(0)}$ is not covered by Theorem \ref{thm-existence}.
In this case, the speed $c_1$ of the generalized transition wave $v$ satisfies
$c_1(t)=2\sqrt{f'(0)}+\xi(t)=\frac{\mu (t)}{\lambda_c} +\lambda_c$.
Thus, in this example, it is possible to improve Theorem \ref{thm-existence}
part 1) to the case $\ul c= 2\sqrt{\ul \mu}$.
\end{ex}

\begin{ex}\label{E:2speeds}
{\em Speeds which cannot be written as
$c(t)=\frac{\mu(t)}{\kappa}+\kappa$.}\\ 
Consider a smooth positive function $\mu=
\mu (t)$ such that $\mu (t)\to \mu_1$ as $t\to -\infty$ and $\mu (t) \to \mu_2$
as $t\to +\infty$, with 
$\mu_1>0$ and $\mu_2>0$. Let $f(t,u)= \mu (t) u (1-u)$. Then it has been proved
by Berestycki and Hamel (in a more general framework, see
\cite{BerestyckiHamelgtw2})
that if $\mu_1<\mu_2$, then for all $c_1\in [2\sqrt{\mu_1},+\infty)$, there
exists a generalized transition wave of equation (\ref{eqprinc}) with speed $c$
such that $c(t)\to c_1$ as $t\to -\infty$ and 
$\frac{1}{t}\int_0^t c(s)ds\to c_2$ as $t\to +\infty$, where $c_2
=\frac{\mu_2}{\kappa_1}+\kappa_1$ and $\kappa_1$ is the smallest root of 
$\kappa_1^2-\kappa_1 c_1 +\mu_1$. This result can be deduced from Theorem
\ref{thm-existence} when $c_1>2\sqrt{\mu_1}$, which even gives a stronger result
since we get $c(t)= \frac{\mu (t)}{\kappa_1}+\kappa_1$ for all $t\in\R$. 

When $\mu_1>\mu_2$, then Berestycki and Hamel obtained a different result. Namely,
they prove that for all $c_1\in [2\sqrt{\mu_1},+\infty)$, 
if $\kappa_1\geq \sqrt{\mu_2}$ (which is true in particular when
$c_1=2\sqrt{\mu_1}$), 
there exists 
a generalized transition wave of equation (\ref{eqprinc}) with speed $c$ such
that $c(t)\to c_1$ as $t\to -\infty$ and 
$\frac{1}{t}\int_0^t c(s)ds\to 2\sqrt{\mu_2}$ as $t\to +\infty$.
In this case the speed $c$ cannot be put in the form $c(t)= \frac{\mu
(t)}{\kappa}+\kappa$ for some $\kappa>0$. Hence, the class of speeds we
construct in Theorem \ref{thm-existence}
is not exhaustive. Moreover, in this example, this class of speeds misses the
most important generalized transition waves: the one which travels with speed 
$2\sqrt{\mu_1}$ when $t\to -\infty$ and $2\sqrt{\mu_2}$ when $t\to +\infty$. As
these two speeds are minimal near $t=\pm\infty$, one can expect this generalized
transition wave to 
be attractive in some sense, as in homogeneous media (see \cite{KPP}).
\end{ex}

%%%%%%%%%%%%%%%%%%%%%%%%%%%%%%%%%%%%%%%%%%%%%%%%%%%%%%%%%%%%%%%%%%%%%%%%%%%%%%%%

%%%%%%%%%%%%%%%%%%%%%%%%%%%%%%%%%%%%%%%%%%%%%%%%%%%%%%%%%%%%%%%%%%%%%%%%%%%%%%%%

\section{Proof of the results}

As we said at the beginning of Section \ref{sec:hyp}, in this paper the terms
(strong)
sub and supersolution refer to functions in $W^{2,1}_{N+1,loc}$, satisfying the
differential inequalities a.e. 
We say that a function
is a {\em generalized subsolution} (resp.~{\em supersolution}) if it is the
supremum (resp.~infimum) of a finite number of subsolutions
(resp.~supersolutions).

%-------------------------------------------------------------------

\subsection{Properties of the least mean}

We first give an equivalent formulation of the least mean
(see Definition \ref{def-mean}). 

\begin{prop} \label{reformulation}
If $g\in L^\infty(\R)$ then its least mean $\ul g$ satisfies $$\ul
g=\lim_{T\to+\infty}\essinf_{t\in\R}\frac{1}{T}\int_t^{t+T}g(s)\,ds.$$ 
In particular, if $g$ admits a mean value 
$\langle g\rangle$, defined by
\eq{mv}, then $\ul g=\langle g\rangle$.
\end{prop}

\begin{proof}
For $T>0$, define the following function:
$$F(T):=\inf_{t\in\R}\int_t^{t+T}g(s)\,ds.$$
We have that
$$\ul g=\sup_{T>0}\frac{F(T)}{T}\geq\limsup_{T\to+\infty}\frac{F(T)}{T}.$$
Therefore, to prove the statement we only need to show that
 $\liminf_{T\to+\infty}F(T)/T\geq\ul g$.
 For any $\e>0$, let $T_\e>0$ be such that $F(T_\e)/T_\e\geq\ul g-\e$.
 We use the notation $\lfloor x \rfloor$ to indicate the floor of the real
 number $x$ (that is, the greatest integer $n\leq x$) and we compute
 $$\fa T>0,\quad F(T)=\inf_{t\in\R}\left(\int_t^{t+\left\lfloor\frac{T}{T_\e}
 \right\rfloor T_\e} g(s)\,ds+\int_{t+\left\lfloor\frac{T}{T_\e}
 \right\rfloor T_\e}^{t+T} g(s)\,ds\right)\geq
 \left\lfloor\frac{T}{T_\e}\right\rfloor F(T_\e)-\|g\|_{L^\infty(\R)}T_\e.$$
 As a consequence,
 $$\liminf_{T\to+\infty}\frac{F(T)}{T}\geq
 \lim_{T\to+\infty}\left\lfloor\frac{T}{T_\e}\right\rfloor
 \frac{T_\e}{T}\frac{F(T_\e)}{T_\e}=\frac{F(T_\e)}{T_\e}\geq\ul g-\e.$$
 The proof is thereby achieved due to the arbitrariness of $\e$.
 \end{proof}

We now derive another characterization of the least mean. 
This is the property underlying the fact that the existence of
generalized transition waves is expressed
in terms of the least mean of their
speeds.

\begin{lem}\label{lem:lm>0}
Let $B\in L^\infty(\R)$. Then 
$$\ul B=\sup_{A\in W^{1,\infty}(\R)}\essinf_{t\in\R}(A'+B)(t).$$
\end{lem}

\begin{proof} If $B$ is a periodic function, then
$g(t):= \underline{B}-
B(t)$ is periodic with zero mean. Thus $A(t):=\int_0^t g(s)ds$ is bounded and
satisfies $A'+B\equiv\underline{B}$. This shows that 
\Fi{lm<}
\ul B\leq\sup_{A\in W^{1,\infty}(\R)}\essinf_{t\in\R}(A'+B)(t)
\Ff
in this simple case. We will now
generalize this construction in order to handle general functions B.

Fix $m<\ul B$. By definition, there exists $T>0$ such that
$$\inf_{t\in\R}\frac1T\int_t^{t+T}B(s)\, ds>m.$$
We define
$$ \forall k\in\Z,\ \text{for a.e.~} t\in
[(k-1)T,kT),\quad
g(t):=-B(t)+\beta_k,\qquad\text{where }\
\beta_k:=\frac{1}{T}\int_{(k-1)T}^{kT}
B(s)\, ds.$$
Then we set $A(t):=\int_0^t g(s)ds$.
It follows that $A'+B\geq m$ and, since $\int_{(k-1)T}^{kT}
g(s)\, ds=0$, that 
$$\|A\|_{L^\infty(\R)}\leq\|g\|_{L^\infty(\R)} T\leq2T\|B\|_{L^\infty(\R)}.$$
Therefore \eq{lm<} holds due to the arbitrariness of $m<\ul B$.

Consider now a function $A\in W^{1,\infty}(\R)$.
Owing to Proposition \ref{reformulation} we derive
\[\begin{split}
\ul B &=\lim_{T\to+\infty}\inf_{t\in\R}\frac{1}{T}\int_t^{t+T}B(s)\,ds
\geq\essinf_\R(A'+B)+
\lim_{T\to+\infty}\inf_{t\in\R}\frac{1}{T}\int_t^{t+T}(-A'(s))\,ds\\
&=\essinf_\R(A'+B)+\lim_{T\to+\infty}\inf_{t\in\R}\frac{1}{T}(A(t)-A(t+T))=
\essinf_\R(A'+B).
\end{split}\]
This concludes the proof.
\end{proof}

\begin{rmq}\label{rmq:lm>0}
 In the proof of Lemma \ref{lem:lm>0} we have shown the following fact: if
 $\eta\in\N\cup\{+\infty\}$, $T>0$ are such that
 $$m:=\inf_{\su{k\in\N}{k\leq\eta}}\frac1T\int_{(k-1)T}^{kT}
 B(s)\, ds>0,$$
 then there exists $A\in W^{1,\infty}((0,\eta T))$ satisfying
 $$\essinf_{[0,\eta T)}(A'+B)=m,\qquad 
 \|A\|_{L^\infty((0,\eta T))}\leq2T
 \|B\|_{L^\infty((0,\eta T))}.$$ 
\end{rmq}

%---------------------------------------------------------------

\subsection{Construction of the generalized transition waves when $\ul c>2\sqrt{\ul \mu}$.}

In order to construct generalized transition waves, we will use appropriate sub
and supersolutions. 
The particular form of the speeds (\ref{formc})
will naturally emerge from constraints on the exponential supersolution.

\begin{prop}\label{pro:subsuper}
Under the assumptions of Theorem \ref{thm-existence}, for all $\gamma>
2\sqrt{\ul\mu}$, there exists a function $c\in L^\infty(\R)$, with 
$\ul c=\gamma$, such that (\ref{def-tf}) admits some uniformly continuous
generalized sub and supersolutions
$\ul\phi(z,t)$, $\ol\phi(z)$ satisfying

$$0\leq\ul\phi<\ol\phi\leq1,\qquad
\ol\phi(+\infty)=0,\qquad\ol\phi(-\infty,t)=1\  \hbox{
uniformly in } t\in\R,$$
$$\exists\xi\in\R,\ \ \inf_{t\in\R}\ul\phi(\xi,t)>0,\qquad
\forall z\in\R,\ \ \inf_{t\in\R}(\ol\phi-\ul\phi)(z,t)>0.$$
$$\ol\phi\ \text{ is nonincreasing in }\R,\qquad
\fa\tau>0,\ \
\lim_{z\to+\infty}\frac{\ol\phi(z+\tau)}{\ul\phi(z,t)}<1
\  \hbox{ uniformly in } t\in\R.$$
\end{prop}

\begin{proof}
Fix $\gamma>2\sqrt{\ul{\mu}}$. We choose $c$ in such a way that the
linearized equation
around $0$ associated with (\ref{def-tf}) admits an exponential solution of the
type $\psi(z)=e^{-\kappa z}$, for some $\kappa>0$. Namely,
$$0=\partial_t \psi-\partial_{zz}\psi -c(t)\partial_z
\psi-\mu(t)\psi=
[-\kappa^2+c(t)\kappa-\mu(t)]\psi,\quad \text{for a.e.~}t\in\R.$$
It follows that $c\equiv\kappa+\kappa^{-1}\mu$. Imposing $\ul c=\gamma$
yields
$$\gamma=\lim_{T\to +\infty}\inf_{t\in\R}
\frac{1}{T}\int_t^{t+T}
[\kappa+\kappa^{-1}\mu(s)]ds=\kappa+\kappa^{-1}\ul\mu.$$
Since $\gamma>2\sqrt{\ul{\mu}}$, the equation
$\kappa^2-\gamma\kappa+\ul\mu=0$
has two positive solutions. We take the smallest one,
that is,
$$\kappa=\frac{\gamma-\sqrt{\gamma^2-4\ul\mu}}2.$$
Extending $f(t,\.)$ linearly outside $[0,1]$, we can assume that $\psi$ is a
global supersolution. We then set $\ol\phi(z):=\min(\psi(z),1)$.

Let $C$, $\gamma$ be the constants in (\ref{hyp-lower}).
Our aim is to find a function $A\in
W^{1,\infty}(\R)$ and a constant $h>\kappa$ such that the
function $\vp$ defined by $\vp(z):=\psi(z)-e^{A(t)-hz}$
satisfies
\begin{equation}\label{ulphi}
\partial_t\vp-\partial_{zz}\vp-c(t)\partial_z\vp
\leq\mu(t)\vp-C\vp^{1+\gamma},\quad\text{for a.e.~} z>0,\ t\in\R.
\end{equation}
By direct computation we see that
$$\partial_t \vp-\partial_{zz}\vp-c(t)\partial_z\vp-\mu(t)\vp
=[-A'(t)+h^2-c(t)h+\mu(t)]e^{A(t)-hz}.$$
Hence (\ref{ulphi}) holds if and only if
$$\text{for a.e.~} z>0,\ t\in\R,\quad
A'(t)+B(t)\geq C\vp^{1+\gamma}e^{hz-A(t)},\qquad\text{where }
B(t):=-h^2+c(t)h-\mu(t).$$
Let $\kappa< h<(1+\gamma)\kappa$. Since
$$\forall z>0,\ t\in\R,\quad
\vp^{1+\gamma}e^{hz-A(t)}\leq
e^{[h-(1+\gamma)\kappa]z-A(t)}
\leq e^{-A(t)},$$
if $\essinf_\R(A'+B)>0$ then the desired inequality follows by adding
a large constant to $A$. 
Owing to Lemma \ref{lem:lm>0}, this condition is fulfilled by a suitable
function $A\in W^{1,\infty}(\R)$ as soon as $\ul B>0$.
Let us compute
\[\begin{split}
\ul B &=
\lim_{T\to +\infty}\inf_{t\in\R}
\frac{1}{T}\int_t^{t+T}h\left[
\kappa-h+\mu(s)\left(\frac1\kappa-\frac1h\right)\right ] ds \\
& =h\left[
\kappa-h+\ul\mu\left(\frac1\kappa-\frac1h\right)\right ]\\
&= -h^2+\gamma h-\ul\mu.
\end{split}\]
Since $\kappa$ is the smallest root of $-x^2+\gamma x-\ul\mu=0$, we can choose
$h\in(\kappa,(1+\gamma)\kappa)$ in such a way that
$\ul B>0$. Therefore, there exists $A\in
W^{1,\infty}(\R)$ such that (\ref{ulphi}) holds. Up to adding a suitable
constant $\alpha<0$ to $A$, it is
not
restrictive to assume that $\vp$ is less than the constant $\delta$ in
(\ref{hyp-lower}). Hence, $\vp$ is a subsolution of (\ref{def-tf}) in
$(0,+\infty)\times\R$.
Since $\vp\geq0$ if and only if $z\geq(h-\kappa)^{-1}A(t)$,
$\alpha$ can be chosen in such a way that $\vp\leq0$ for $z\leq0$. 
% Then, defining $f(t,u):=\mu(t)u$ for $u\leq0$, we obtain
% $$\partial_t\vp-\partial_{zz}\vp-c(t)\partial_z\vp-f(t,\vp)
% =-(A'(t)+B(t))e^{A(t)-hz}\leq0\quad\text{in }(-\infty,0)\times\R.$$
Whence, due to the arbitrariness of the extension of $f(t,u)$ for $u<0$, it
follows that $\ul\phi(z,t):=\max(\vp(z,t),0)$ is a generalized
subsolution. 
Finally, for $\tau>0$, it holds that
$$\lim_{z\to+\infty}\frac{\ol\phi(z+\tau)}{\ul\phi(z,t)}=
\lim_{z\to+\infty}\frac{e^{-\kappa\tau}}{1-e^{A(t)-(h-\kappa)z}}=e^{-\kappa\tau}
.$$
\end{proof}

\begin{proof}[Proof of Theorem \ref{thm-existence} part 1)]
Let $c,\ \ul\phi,\ \ol\phi$ be the functions given by Proposition
\ref{pro:subsuper}, with $\ul c=\gamma$. For $n\in\N$, consider the solution
$\phi_n$ of the problem
\Fi{-n}
\left\{\begin{array}{ll}
\partial_t \phi-\partial_{zz}\phi -c(t)\partial_z \phi = f(t,\phi),& 
z\in\R,\ t>-n, \\ 
\phi(z,-n)=\ol\phi(z), & z\in\R. 
\end{array} \right.
\Ff
The comparison principle implies that $\ul\phi\leq\phi_n\leq\ol\phi$ and, since
$\ol\phi$ is nonincreasing, that $\phi_n(\.,t)$ is nonincreasing too.
Owing to the parabolic estimates and the embedding theorem, using
a diagonal extraction method we can
find a subsequence of $(\phi_n)_{n\in\N}$ converging  weakly in $W^{2,1}_p(K)$
and 
strongly in $L^\infty(K)$, for any compact $K\subset\R\times\R$ and any
$p<\infty$, to a solution $\phi$ of 
$$\partial_t \phi-\partial_{zz}\phi -c(t)\partial_z \phi = f(t,\phi),\quad
z\in\R,\ t\in\R.$$ 
The function $\phi$ is nonincreasing in $z$ and
satisfies $\ul\phi\leq\phi\leq\ol\phi$. Applying the parabolic strong maximum
principle to $\phi(z-z_0,t)-\phi(z,t)$, for every $z_0>0$, we find that $\phi$
is decreasing in~$z$.

It remains to prove that $\phi(-\infty,t)=1$ uniformly with respect to $t\in\R$.
Set
$$\theta:=\lim_{z\to-\infty}\inf_{t\in\R}\phi(z,t).$$
Our aim is to show that $\theta=1$.
Let $(t_n)_{n\in\N}$ be such that $\lim_{n\to\infty}\phi(-n,t_n)=\theta$. 
We would like to pass to the limit in
the sequence of equations satisfied by the $\phi(\.-n,\.+t_n)$, but this is not
possible due to the presence of the drift term. To
overcome this difficulty we come back to the fixed coordinate system by
considering
the functions $(v_n)_{n\in\N}$ defined by
$$v_n(z,t):=\phi (z-n-\int_{t_n}^{t_n+t}c(s)ds,t+t_n).$$
These functions are solutions of
$$\partial_t v_n-\partial_{zz}v_n= f(t+t_n,v_n),\quad
z\in\R,\ t\in\R,$$ 
and satisfy $\lim_{n\to\infty}v_n(0,0)=\theta$ and 
$\liminf_{n\to\infty}v_n(z,t)\geq\theta$ locally uniformly in
$(z,t)\in\R\times\R$.
The same diagonal extraction method as before 
shows that $(v_n)_{n\in\N}$ converges (up to subsequences) weakly in
$W^{2,1}_{p,loc}$ and strongly in $L^\infty_{loc}$ to some function $v$ 
satisfying
$$\partial_t v-\partial_{zz}v=g(z,t)\geq0,\quad\text{for a.e.~}z\in\R,\
t\in\R,$$
where $g(z,t)$ is the weak limit in $L^p_{loc}(\R\times\R)$ of (a subsequence
of) $f(t+t_n,v_n)$.
Moreover, $v$ attains its minimum value $\theta$ at $(0,0)$.
As a consequence, the strong maximum principle yields
$v=\theta$ in $\R\times(-\infty,0]$. In particular, $g=0$ a.e.~in
$\R\times(-\infty,0)$.
Using the Lipschitz continuity of $f(t,\.)$,
we then derive
$$\text{for a.e.~}(z,t)\in\R\times(-\infty,0),\quad
0=g(z,t)\geq\essinf_{s\in\R}f(s,\theta).$$
Therefore, hypothesis \ref{hyp-pos} yields 
$\theta=0$ or $1$.
The proof is then concluded by noticing that
$$\theta=\lim_{z\to-\infty}\inf_{t\in\R}\phi(z,t)
\geq\inf_{t\in\R}\phi(\xi,
t)\geq\inf_{t\in\R}
\ul\phi(\xi,t)>0$$
($\xi$ being the constant in Proposition \ref{pro:subsuper}).
\end{proof}

%--------------------------------------------------------------------------

\subsection{Non-existence of generalized transition waves when $\ul c <
2\sqrt{\ul \mu}$.}\label{sec:nonE}

This section is dedicated to the proof of the lower bound for the least mean of
admissible speeds - Theorem \ref{thm-existence} part 2.
This is achieved by comparing the generalized transition waves with some
subsolutions whose level-sets propagate at speeds less than $2\sqrt{\ul
\mu}$. The construction of the
subsolution is based on an auxiliary result which
is quoted from \cite{BHRossi} and reclaimed in the Appendix here.

\begin{lem}\label{lem:0nT}
 Let $g\in L^\infty(\R)$, $\eta\in\N\cup\{+\infty\}$ and $T>0$ be such that
 $$\gamma_*:=\inf_{\su{k\in\N}{k\leq\eta}}\,
 2\sqrt{\frac1T\,\int_{(k-1)T}^{kT}g(s)\,
 ds}>0.$$
 Then for all $0<\gamma<\gamma_*$, there exists a uniformly continuous
 subsolution $\ul v$ of 
 \Fi{0etaT}
 \partial_t v-\Delta v=g(t)v, \quad x\in\R^N,\ t\in(0,\eta T),
 \Ff
 such that 
 $$0\leq\ul v\leq1,\qquad \ul v(x,0)=0 \text{ for }|x|\geq R,\qquad
 \inf_{\su{0\leq t<\eta T}{|x|\leq\gamma t}}\ul v(x,t)\geq C,$$
 where $R$, $C$ only depend on $T$, $\gamma_*-\gamma$, $N$,
 $\|g\|_{L^\infty((0,\eta T))}$ and not on $n$.
\end{lem}

\begin{proof}
Fix $\gamma\in(0,\gamma^*)$. By Lemma \ref{lem:h} in the Appendix,
there exist $h\in C^{2,\alpha}(\R)$ and $r>0$, both
depending on $\gamma$ and $\gamma_*-\gamma$, satisfying
$$h=0\ \text{in }(-\infty,0],\qquad
h'>0\ \text{in }(0,r),\qquad 
h=1\text{ in }[r,+\infty),$$
$$Q\leq \gamma+1,\quad
4C-Q^2\geq\frac12(\gamma_*-\gamma)^2
\qquad\Rightarrow\qquad-h''+Qh'-Ch<0\quad\text{in }(0,r).$$
Note that $r$, $h$ actually depend on
$\gamma_*-\gamma$ and $\|g\|_\infty$, because $\gamma<\gamma_*\leq\|g\|_\infty$.
We set
$$\ul v(x,t):=e^{A(t)}h(R-|x|+\gamma t),$$
where $R>r$ and $A:\R\to\R$ will be chosen later. This function solves
$$\partial_t \ul v-\Delta\ul v-g(t)\ul v
=\left[-h''(\rho)+\left(\gamma+\frac{N-1}{|x|}\right)h'(\rho)-C(t)h(\rho)\right]
e^{A(t)},\quad x\in\R^N,\ t\in(0,\eta T).$$
with $\rho=R-|x|+\gamma t$ and $C(t)=g(t)-A'(t)$.
Since $h'$ is nonnegative and vanishes in $[r,+\infty)$, it
follows that
$$\partial_t \ul v-\Delta\ul v-g(t)\ul v
\leq[-h''(\rho)+Qh'(\rho)-C(t)h(\rho)]e^{A(t)},\quad \text{for a.e.~} x\in\R^N,\
t\in(0,\eta T),$$
with $Q=\gamma+\frac{N-1}{R-r}$. We write
$$4C(t)-Q^2=B(t)-4A'(t),\quad\text{where }
B(t)=4 g(t)-\left(\gamma+\frac{N-1}{R-r}\right)^2.$$
and we compute
$$m:=\inf_{\su{k\in\N}{k\leq\eta}}\frac1T
\int_{(k-1)T}^{kT}B(s)\,
ds=\gamma_*^2-\left(\gamma+\frac{N-1}{R-r}\right)^2.$$
Hence, since $\gamma_*^2-\gamma^2\geq(\gamma_*-\gamma)^2$, it is possible to
choose $R$, depending on $N$, $r$, $\|g\|_\infty$ and $\gamma_*-\gamma$, such
that $m\geq\frac12(\gamma_*-\gamma)^2$.
By Remark \ref{rmq:lm>0} 
there exists a function $A\in W^{1,\infty}(\R)$ such that
$$\min_{[0,\eta T)}(B-4A')=m,\qquad \|A\|_{L^\infty((0,\eta T))}\leq\frac T2
\|B\|_{L^\infty((0,\eta T))}\leq4T\|g\|_{L^\infty((0,\eta T))}.$$
Consequently, $4C-Q^2\geq\frac12(\gamma_*-\gamma)^2$ a.e.~in $(0,\eta T)$ and,
up to increasing $R$ if need be, $Q\leq\gamma+1$.
Therefore, $\ul v$ is a subsolution of \eq{0etaT}. This concludes the proof.
\end{proof}

\begin{proof}[Proof of Theorem \ref{thm-existence} part 2)]
Let $u$ be a generalized transition wave with speed $c$. 
Definition \ref{def-gtf} yields
\Fi{Le}
\lim_{L\to+\infty}\,\inf_{\su{x\.e<-L}{t\in\R}}u(x+e\int_0^tc(s)ds,t)=1,
\qquad
\lim_{L\to+\infty}\,\sup_{\su{x\.e>L}{t\in\R}}u(x+e\int_0^tc(s)ds,t)=0.
\Ff
By the definition of least mean, for all $\e>0$, there exists $T>0$ such
that 
$$\fa T'\geq T,\quad
\frac1{T'}\inf_{t\in\R}\int_t^{t+T'}c(s)\,ds<\ul c+\e,\qquad
\fa t\in\R,\quad\frac1T\int_t^{t+T}\mu(s)\,ds>\ul\mu-\e.$$
For $n\in\N$, let $t_n$ be such that
$$\frac1{nT}\int_{t_n}^{t_n+nT}c(s)\,ds<\ul c+2\e.$$
Taking $\e$ small enough in such a way that $2\sqrt{\ul\mu-2\e}>\e$, we find
that 
$$\gamma_*^n:=\inf_{k\in\{1,\dots,n\}}2\sqrt{\frac1T\,\int_{(k-1)T}^{kT}
(\mu(s+t_n)-\e)\,
ds}>2\sqrt{\ul\mu-2\e}>\e.$$
Let $(\ul v_n)_{n\in\N}$ be the functions obtained applying Lemma
\ref{lem:0nT} with $g(t)=\mu(t+t_n)-\e$
and $\gamma=\gamma^n:=\gamma_*^n-\e$, and
let $R$, $C$ be the associated constants, which are independent of $n$.
By the regularity hypothesis on $f$, there exists $\sigma\in(0,1)$ such that
$f(t,w)\geq(\mu(t)-\e)w$ for $w\in[0,\sigma]$. As a consequence, the functions
$\ul u_n(x,t):=\sigma \ul v_n(x,t-t_n)$ satisfy
$$\partial_t \ul u_n-\Delta \ul u_n\leq f(t,\ul u_n), \quad \text{for a.e.~}
x\in\R^N,\
t\in(t_n,t_n+nT).$$
By \eq{Le}, for $L$ large enough we have that
$$\fa t\in\R,\quad\inf_{|x|<R}u(x-Le+e\int_0^tc(s)ds,t)\geq\sigma.$$
Thus, up to replacing $u(x,t)$ with $u(x-Le,t)$,
it is not restrictive to assume that
$$\fa n\in\N,\ x\in\R^N,\quad u(x+e\int_0^{t_n}c(s)ds,t_n)\geq \ul u_n(x,t_n).$$
The comparison principle then yields
$$\fa n\in\N,\ x\in\R^N,\ t\in(t_n,t_n+nT),\quad
u(x+e\int_0^{t_n}c(s)ds,t)\geq \ul u_n(x,t).$$
Therefore, 
$$\liminf_{n\to\infty}u(\gamma^n nT e+e\int_0^{t_n}c(s)ds,t_n+nT)\geq
\liminf_{n\to\infty}\ul u_n(\gamma^n nT,t_n+nT)\geq\sigma C,$$
whence, owing to \eq{Le}, we deduce that
\[\begin{split}
+\infty &>\limsup_{n\to\infty}\left(\gamma^n nT +\int_0^{t_n}c(s)ds
-\int_0^{t_n+nT}c(s)ds\right)\\
&=\limsup_{n\to\infty}\left(\gamma^n nT -\int_{t_n}^{t_n+nT}c(s)ds\right)\\
&\geq \limsup_{n\to\infty}\left(2\sqrt{\ul\mu-2\e}-\e-\ul c-2\e\right)nT.
\end{split}\]
That is, $\ul c\geq2\sqrt{\ul\mu-2\e}-3\e$.
Since $\e$ can be chosen arbitrarily small, we eventually infer that $\ul
c\geq2\sqrt{\ul\mu}$.
\end{proof}

\begin{rmq}
The same arguments as in the above proof yield the non-existence of fronts
$(\phi,c)$ such that $\ul c_{\pm}<2\sqrt{\ul\mu_\pm}$, where, 
for a given function $g\in L^\infty(\R)$,
$$\underline{g}_+:=\sup_{T>0}\,\inf_{t>0}
\frac{1}{T}\int_t^{t+T}g(s)\,ds,\qquad
\underline{g}_-:=\sup_{T>0}\,\inf_{t<0}
\frac{1}{T}\int_{t-T}^T g(s)\,ds$$
\end{rmq}

%%%%%%%%%%%%%%%%%%%%%%%%%%%%%%%%%%%%%%%%%%%%%%%%%%%%%%%%%%%%%%%%%%%%%%%%%%%%%%%%

\section{The random stationary ergodic case}

To start with, we show that the temporal least mean of a stationary ergodic
function
is almost surely independent of $\omega$.

\begin{prop}\label{pro:lmergodic}
For a given bounded
measurable function $\tilde{g}:\O\to\R$, the mapping 
$$\omega\,\mapsto\,\sup_{T>0}\,\inf_{t\in\R}
\frac{1}{T}\int_t^{t+T}\tilde{g}(\pi_s\omega)\,ds$$
is constant in a set of probability measure $1$.
We call this constant value the {\em least mean} of the random stationary ergodic function $g$ defined by $g(t,\omega):= \tilde{g}(\pi_t \omega)$
and we denote it by $\ul g$.
\end{prop}

\begin{proof}
The result follows immediately from the ergodicity of the process
$(\pi_t)_{t\in\R}$. Indeed, setting 
$$G(\omega):=\sup_{T>0}\,\inf_{t\in\R}\frac{1}{T}\int_t^{t+T}\tilde{g}(\pi_s\omega)\,ds,$$
for given $\e>0$, there exists a set $\mc{A}_\e\in\mc{F}$, with
$\P(\mc{A}_\e)>0$, such that $G(\omega)<\essinf_\O G+\e$ for
$\omega\in\mc{A}_\e$.
It is easily seen that $\mc{A}_\e$ is invariant under the action of
$(\pi_t)_{t\in\R}$, and then $\P(\mc{A}_\e)=1$. Owing to the arbitrariness of
$\e$, we infer that $G$ is almost surely equal to $\essinf_\O G$.
\end{proof}

The proof of Theorem \ref{thm-existencerandom} relies on a general uniqueness
result for the profile of  generalized transition waves that share the 
same behavior at infinity. This, in turn, is derived using the following strong
maximum principle-type property.

\begin{lem}\label{lem:I}
Let $c\in L^\infty(\R)$ and assume that $f$ satisfies the regularity
conditions
in Hypothesis \ref{hyp:KPP}.
Let $I$ be an open interval and $\vp,\psi$ be respectively a generalized sub
and supersolution of 
$$\partial_t\phi-\partial_{zz}\phi-c(t)\partial_z\phi
=f(t,\phi),\quad z\in I,\ t\in\R,$$ 
which are uniformly continuous and
satisfy $0\leq\vp\leq\psi\leq1$ in $I\times\R$. Then, the set
$$\{z\in I\ : \ \inf_{t\in\R}(\psi-\vp)(z,t)=0\}$$
is either empty or coincides with $I$.
\end{lem}

\begin{proof}
Clearly, it is sufficient to prove the result for strong sub and
supersolutions. We achieve this by showing that the set
$$J:=\{z\in I\ : \ \inf_{t\in\R}(\psi-\vp)(z,t)=0\}$$
is open and closed in the topology of $I$.
That it is closed follows immediately from the uniform continuity of $\vp$ and
$\psi$.

Let us show that it is open. Suppose that there exists $z_0\in J$. There
is a sequence $(t_n)_{n\in\N}$ such that $(\psi-\vp)(z_0,t_n)$ tends to $0$ as
$n$ goes to infinity. For $n\in\N$, define $\Phi_n(z,t):=(\psi-\vp)(z,t+t_n)$.
These functions satisfy
$$\partial_t\Phi_n-\partial_{zz}\Phi_n-c(t+t_n)\partial_z\Phi_n
-\zeta(z,t+t_n)\Phi_n\geq0\quad \text{for a.e.~}z\in I,\ t\in\R,$$
where
$$\zeta(z,t):=\frac{f(t,\psi)-f(t,\vp)}{\psi-\vp}$$
belongs to $L^\infty(I\times\R)$
due to the Lipschitz-continuity of $f$. Let $\delta>0$ be such that
$[z_0-\delta,z_0+\delta]\subset I$. We now make use of the parabolic weak
Harnack
inequality (see e.g. Theorem 7.37 in \cite{Lie}). It provides two constants
$p,C>0$
such that
$$\fa n\in\N,\quad
\|\Phi_n\|_{L^p((z_0-\delta,z_0+\delta)\times(-2,-1))}
\leq C\inf_{(z_0-\delta,z_0+\delta)\times(-1,0)}\Phi_n
\leq C\Phi_n(z_0,0).$$
Whence $(\Phi_n)_{n\in\N}$ converges to $0$ in
$L^p((z_0-\delta,z_0+\delta)\times(-2,-1))$. By the Arzela-Ascoli
theorem we then infer that, up to subsequences, $(\Phi_n)_{n\in\N}$ converges to
$0$ uniformly in $(z_0-\delta,z_0+\delta)\times(-2,-1)$. This means that
$(z_0-\delta,z_0+\delta)\subset J$.
\end{proof}

\begin{prop}\label{pro:!}
Assume that $c\in L^\infty(\R)$ and that $f$ satisfies (\ref{hyp-vois1}). Let
$\vp,\ \psi$ be a subsolution and a positive supersolution of (\ref{def-tf})
which are uniformly continuous and satisfy
$$0\leq\vp,\psi\leq1,\qquad
\psi(\.,t) \text{ is nonincreasing},\qquad\fa\tau>0,\quad
\lim_{R\to+\infty}\sup_{\su {z>R}{t\in\R}}\frac{\vp(z,t)}{\psi(z-\tau,t)}<1.$$
Then $\vp\leq\psi$ in $\R\times\R$.
\end{prop}

\begin{proof}
Since $\vp,\ \psi$ are uniformly continuous and satisfy $\vp(+\infty,t)=0$,
$\psi(-\infty,t)=1$ uniformly in $t$, applying Lemma \ref{lem:I}, first with
$\psi\equiv1$ and then with $\vp\equiv0$, we derive 
\Fi{supinf}
\fa r\in\R,\quad\sup_{(r,+\infty)\times\R}\vp<1,\quad
\inf_{(-\infty,r)\times\R}\psi>0.
\Ff
Let $\delta\in(0,1)$ be the constant in (\ref{hyp-vois1}).
By hypothesis, there exists $\rho\in\R$ such that $\psi>1-\delta$ in
$(-\infty,\rho)\times\R$. Let $\chi:\R\to[0,1]$ be a smooth function satisfying
$\chi=1$ in $(-\infty,\rho]$, $\chi=0$ in $[\rho+1,+\infty)$. Define the family
of functions $(\psi^{\e,\tau})_{\e,\tau\geq0}$ by setting
$$\psi^{\e,\tau}(z,t):=[1+\e\chi(z)]\psi(z-\tau,t).$$
Since $\lim_{\e,\tau\to0^+}\psi^{\e,\tau}\equiv\psi$, the statement is proved
if we show that $\psi^{\e,\tau}\geq\vp$ for all $\e,\tau>0$.
The $\psi^{\e,\tau}$ are nondecreasing with respect to both $\e$ and $\tau$.
Moreover, there
exists $z_0\in\R$ such that $\psi^{0,1}>\vp$ in $(z_0,+\infty)\times\R$. On the
other hand, for all $\e>0$, there exists $\tau\geq0$ such that
$\psi^{\e,\tau}>1\geq\vp$ in $(-\infty,z_0]\times\R$. Consequently, for all
$\e>0$, we have that $\psi^{\e,\tau}\geq\vp$ for $\tau$ large enough. Define
$$\fa\e>0,\quad\tau(\e):=\min\{\tau\geq0\ : \ \psi^{\e,\tau}\geq\vp\}.$$
The function $\e\mapsto\tau(\e)$ is nonincreasing and it holds that
$\psi^{\e,\tau(\e)}\geq\vp$. We argue by contradiction, assuming that there
exists $\t\e>0$ such that $\tau(\t\e)>0$. 
By hypothesis, we have that
\Fi{psi>vp}
\exists\, h>1,\ R\in\R,\quad
\psi^{0,\tau(\t\e)/2}\geq h\vp\quad\text{in }(R,+\infty)\times\R.
\Ff
Fix $\e\in(0,\t\e]$. We know that
$\tau(\e)\geq\tau(\t\e)>0$ and, for $\tau\in(0,\tau(\e))$,
$\inf_{\R\times\R}(\psi^{\e,\tau}-\vp)<0$.
Hence, from \eq{psi>vp} it follows that, for $\tau\in[\tau(\t\e)/2,\tau(\e))$, 
$\inf_{(-\infty,R]\times\R}(\psi^{\e,\tau}-\vp)<0$. Thus, by the uniform
continuity of $\psi$ we get
$\inf_{(-\infty,R]\times\R}(\psi^{\e,\tau(\e)}-\vp)=0$.
We now use the assumption (\ref{hyp-vois1}). Since
$\psi^{0,\tau(\e)}\geq\psi>1-\delta$ in $(-\infty,\rho)\times\R$, 
for a.e.~$z\in(-\infty,\rho)$, $t\in\R$ we have
$$\partial_t \psi^{\e,\tau(\e)}-\partial_{zz}\psi^{\e,\tau(\e)}
-c(t)\partial_z
\psi^{\e,\tau(\e)} =(1+\e)f(t,\psi^{0,\tau(\e)})\geq
f(t,\psi^{\e,\tau(\e)})$$
(where we have extended $f$ by $0$ in $\R\times(1,+\infty)$).
By hypothesis, we can find $R_\e<\rho$ such that
$\inf_{(-\infty,R_\e]\times\R}\psi^{\e,\tau(\e)}>1\geq\sup\vp$.
Consequently, Lemma \ref{lem:I} yields
$\inf_{(-\infty,\rho-1]\times\R}(\psi^{\e,\tau(\e)}-\vp)>0$.
It follows that $R>\rho-1$ and that
\Fi{psi=vp}
\inf_{(\rho-1,R]\times\R}(\psi^{\e,\tau(\e)}-\vp)=0.
\Ff
In order to pass to the limit $\e\to0^+$ in the above expression, we notice
that, by \eq{supinf}, there exists $\tau_0>0$ such that
$$\inf_{(-\infty,R]\times\R}\psi^{0,\tau_0}>\sup_{[\rho-1,R]\times\R}\vp.$$
As a consequence, \eq{psi=vp} implies that the nonincreasing function $\tau$
is bounded form above by $\tau_0$, and then there exists
$\tau^*:=\lim_{\e\to0^+}\tau(\e)$.
Letting $\e\to0^+$ in the inequality
$\psi^{\e,\tau(\e)}\geq\vp$ and in \eq{psi=vp} yields 
$$\psi^{0,\tau^*}\geq\vp\ \text{
in }\R\times\R,\qquad\inf_{[\rho-1,R]\times\R}(\psi^{0,\tau^*}
-\vp)=0.$$
Thus, using once again Lemma \ref{lem:I} and then the inequality \eq{psi>vp}
we derive
$$0=\inf_{t\in\R}(\psi^{0,\tau^*}-\vp)(R,t)\geq\inf_{t\in\R}(\psi^{0,
\tau(\t\e)/2} -\vp)(R , t)
\geq(1-h^{-1})\inf_{t\in\R}\psi^{0,\tau(\t\e)/2}(R,t).$$
This contradicts \eq{supinf}.
We have shown that $\tau(\e)=0$ for all $\e>0$. That is, $\psi^{\e,\tau}\geq\vp$
for all $\e,\tau>0$. 
\end{proof}

We are now in position to prove Theorem \ref{thm-existencerandom}. 

\begin{proof}[Proof of Theorem \ref{thm-existencerandom}]
First, we fix
$\omega\in\Omega$ such that $\mu(\cdot,\omega)$ admits
$\underline{\mu}$ as a least mean.
By Theorem \ref{thm-existence} there exists a generalized transition wave in
direction $e$ with a speed $c(t,\omega)$ such that $\ul c(\.,\omega)= \gamma$ 
and a profile $\phi(z,t,\omega)$ which is decreasing with respect to $z$.
Moreover, $c(t,\omega)=\kappa+\kappa^{-1}\mu(t,\omega)$, where $\kappa$ is the
unique solution in $(0,\sqrt{\underline{\mu}})$ of $\kappa+\kappa^{-1}
\underline{\mu} =\gamma$.
For $s\in\R$, we set $\phi^s(z,t,\omega):=\phi^s(z,t-s,\pi_s\omega)$. As $f$ and
$c$ are random stationary, the functions $\phi$ and $\phi^s$ satisfy the 
same equation
$$\partial_t \phi-\partial_{zz}\phi
-c(t,\omega)\partial_z\phi = f(t,\omega,\phi),\quad z\in\R, t\in\R.$$
We further know that $\ul\phi\leq\phi\leq\ol\phi$, where
$\underline{\phi}=\underline{\phi}(z,t,\omega)$ and
$\overline{\phi}=\overline{\phi}(z)$ are given by Proposition
\ref{pro:subsuper}. We point out that $\ol\phi(z)\equiv\min(e^{-\kappa z},1)$
does not depend on $\omega$.
For $\tau>0$, we get
$$\lim_{R\to+\infty}\sup_{\su
{z>R}{t\in\R}}\frac{\phi^s(z,t,\omega)}{\phi(z-\tau,t,\omega)}\leq
\lim_{R\to+\infty}\sup_{\su {z>R}{t\in\R}}\frac{\ol\phi(z)}
{\ul\phi(z-\tau,t,\omega)}<1.$$
Hence,
Proposition \ref{pro:!} gives $\phi^s(\cdot,\cdot,\omega)\leq\phi
(\cdot,\cdot,\omega)$. Exchanging the roles of $\phi^s$, $\phi$ we derive
$\phi^s(\cdot,\cdot,\omega)\equiv\phi(\cdot,\cdot,\omega)$ for almost every $\omega$, that is,
$$\fa (z,t,s,\omega)\in\R\times\R\times \R\times \Omega,\quad  \phi
(z,t+s,\omega)=\phi (z,t,\pi_s\omega).$$
Set $\tilde{\phi} (z,\omega):= \phi (z,0,\omega)$ and $\t
c(\omega):=c(0,\omega)$. We see that the function $u$ defined by
$u(x,t,\omega):=
\t\phi (x\cdot e -\int_0^t \tilde{c}(\pi_s\omega)ds,\pi_t\omega)$ is a random
transition wave with random speed $\tilde{c}$ and random
profile $\tilde{\phi}$.
\end{proof}

%%%%%%%%%%%%%%%%%%%%%%%%%%%%%%%%%%%%%%%%%%%%%%%%%%%%%%%%%%%%%%%%%%%%%%%%%%%%

\section{Proof of the spreading properties}

This Section is dedicated to the proof of Proposition \ref{prop-spreading}.
We prove separately properties \eq{spreadinginf} and \eq{spreadingsup}.

\begin{proof}[Proof of \eq{spreadinginf}] 
Let $\gamma<\gamma'<2\sqrt{\ul\mu}$.
Take $0<\e<\ul\mu$ in such a way that
$\gamma'<2\sqrt{\ul\mu-\e}$. By
definition of least mean over $\R_+$ there exists $T>0$ such that 
$$2\sqrt{\inf_{t>0}\frac1T\int_t^{t+T}(\mu(s)-\e)ds}>\gamma'.$$
Hence 
$$\gamma^*:=\inf_{k\in\N}\,2\sqrt{\frac1T\int_{(k-1)T}^{kT}(\mu(s+1)-\e)ds}
>\gamma'.$$
Let $\ul v$ be the subsolution given by Lemma \ref{lem:0nT} with
$g(t)=\mu(t+1)-\e$, $\eta=+\infty$ and $\gamma$ replaced by $\gamma'$.
Since $\ul v(\.,0)$ is compactly supported and $u(\.,1)$ is positive by the
strong maximum principle, it is possible to normalize $\ul v$ in such a way that
$\ul v(\.,0)\leq u(\.,1)$.
Moreover, by further decreasing $\ul v$ if need be, it is not restrictive
to assume that
$$\partial_t\ul  v-\Delta\ul  v=(\mu(t+1)-\e)\ul v\leq f(t+1,\ul v), \quad
\text{for a.e.~} x\in\R^N,\ t>0.$$
Therefore,
$u(x,t+1)\geq \ul v (x,t)$ for $x\in\R^N$, $t>0$ by the comparison principle.
Whence, 
\Fi{m}
m:=\liminf_{t\to +\infty}\inf_{|x|\leq \gamma't} u(x,t+1)
\geq\liminf_{t\to +\infty}\inf_{|x|\leq \gamma' t}\ul v(x,t)>0.
\Ff
If $f$ was smooth with respect to $t$, then the conclusion would follow from Hypothesis \ref{hyp:KPP} (condition
(\ref{hyp-pos}) in particular) through classical arguments. But as we only
assume $f$ to be bounded and measurable in $t$, some more arguments are needed
here.

Set 
$$\theta:=\liminf_{t\to +\infty}\inf_{|x|\leq \gamma t} u(x,t),$$
and let $(x_n)_{n\in\N}$ and $(t_n)_{n\in\N}$ be such that
$$\lim_{n\to\infty}t_n=+\infty,\qquad
\forall n\in\N,\ \ |x_n|\leq\gamma t_n,
\qquad\lim_{n\to\infty}u(x_n,t_n)=\theta.$$
Usual arguments show that, as $n\to\infty$, the functions
$u_n(x,t):=u(x+x_n,t+t_n)$ converge (up to subsequences) weakly in
$W^{2,1}_{p,loc}(\R^N\times\R)$, for any $p<\infty$, and strongly in
$L^\infty_{loc}(\R^N\times\R)$ to a solution $v$ of
$$\partial_t v-\Delta v= g(x,t),\quad x\in\R^N,\ t\in\R,$$
where $g$ is the weak limit in $L^p_{loc}(\R^N\times\R)$ of 
$f(t+t_n,u_n)$. We further see that $v(0,0)=\theta$ and that $v$ is bounded
from below by the constant $m$ in \eq{m}. Let $(\xi_k,\tau_k)_{k\in\N}$ be a
minimizing sequence for $v$, that is,
$$\lim_{k\to\infty}v(\xi_k,\tau_k)=\eta:=\inf_{\R^N\times\R}v.$$
In particular, $0<m\leq\eta\leq\theta$.
The same arguments as before, together with the strong maximum principle, imply
that the $v_k(x,t):=v(x+\xi_k,t+t_k)$ converge (up to subsequences) to
$\eta$ in, say, $L^\infty(B_1\times(-1,0))$.
For $h\in\N\backslash \{0\}$, let $k_h, n_h\in\N$ be such that
$$\|v_{k_h}-\eta\|_{L^\infty(B_1\times(-1,0))}<\frac1h,\qquad
\|u_{n_h}-v\|_{L^\infty(B_1(\xi_{k_h})\times
(\tau_{k_h}-1,\tau_{k_h}))}<\frac1h.$$
Hence, the functions $\t
u_h(x,t):=u(x+x_{n_h}+\xi_{k_h},
t+t_{n_h}+\tau_{k_h})$ satisfy
$\|\t u_h-\eta\|_{L^\infty(B_1\times(-1,0))}<\frac2h$ and then $(\t
u_h)_{h\in\N}$ converges to $\eta$ uniformly in $B_1\times(-1,0)$. On the other
hand, it converges (up to subsequences) to a solution $\t v$ of 
$$\partial_t\t v-\Delta\t  v=\t g(x,t),\quad x\in B_1,\ t\in(-1,0),$$
where $\t g$ is the weak limit in $L^p(B_1\times(-1,0))$ of 
$f(t+t_{n_h}+\tau_{k_h},\t u_h)$. As a consequence,
$$\text{for a.e.~}x\in B_1,\ t\in(-1,0),\quad 0=\t g(x,t)\geq\essinf_{s\in\R}
f(s,\eta).$$
Hypothesis \ref{hyp-pos} then yields $1=\eta\leq\theta$.
\end{proof}

\begin{proof}[Proof of \eq{spreadingsup}]
Let $R>0$ be such that $\supp u_0\subset B_R$. For all $ \kappa>0$ and
$e\in\mathbb{S}^{N-1}$, we define
$$v_{\kappa,e} (x,t):= e^{-\kappa (x\cdot e-R-\kappa t)+\int_0^t \mu (s)ds}. $$
Direct computation shows that the functions $v_{\kappa,e}$ satisfy
$$
\partial_t v_{\kappa,e} -\Delta v_{\kappa,e} -\mu (t)
v_{\kappa,e} =0,\quad x\in\R^N,\ t>0,
\qquad\fa x\in B_R,\quad v_{\kappa,e}(x)>1.$$
Hence, by (\ref{hyp-KPP}), they are supersolutions of \eq{Cauchy} and then
they are greater than $u$ due to the comparison principle.
Let $\sigma>0$, $x\in\R^N$ and $t>0$ be such that
$|x|\geq2\sqrt{t\int_0^t\mu(s)ds}+\sigma t$. Applying the inequality $u(x,t)\leq
v_{\kappa,e}(x,t)$ with $e=\frac x{|x|}$ and $\kappa=\frac{|x|-R}{2t}$ yields
$$u(x,t)\leq\exp\left(-\frac{(|x|-R)^2}{4t}
+\int_0^t\mu(s)ds\right).$$
If in addition $t>R/\sigma$ then $|x|-R\geq2\sqrt{t\int_0^t\mu(s)ds}+\sigma
t-R>0$, whence
$$u(x,t)\leq\exp\left(-\frac{\sigma^2t}{4}-\frac{R^2}{4t}-\sigma\sqrt{t\int_0^t\mu(s)ds}
+R\sqrt{\frac1t\int_0^t\mu(s)ds}+\frac{R\sigma}2\right).$$
Since the right hand side tends to $0$ as $t\to+\infty$, \eq{spreadingsup}
follows.
\end{proof}

%%%%%%%%%%%%%%%%%%%%%%%%%%%%%%%%%%%%%%%%%%%%

\appendix
\section{Appendix}\label{appendix}

The next result, quoted from \cite{BHRossi}, is the key tool used in Section
\ref{sec:nonE} to construct a compactly supported subsolution $\ul v$ to
(\ref{eqprinc}).
For the reader's ease, we include its proof below.
It is slightly simpler than the original one of \cite{BHRossi}.

\begin{lem}[Lemma 8.1 in \cite{BHRossi}]\label{lem:h}
For any given positive numbers $\beta$, $\sigma$, $\theta$, there exist
a function $h\in C^2(\R)$ and a constant $r> 0$ such that
$$h = 0 \text{ in } (-\infty,0],\qquad h' > 0 \text{ in }
(0,r),\qquad h=1 \text{ in }
[r,+\infty),$$
$$-Ah''+Qh'-Ch<0\quad\text{in }(0,r),$$
for all nonnegative constants $A$, $Q$, $C$ satisfying
$$A\leq\beta,\qquad Q\leq\sigma,\qquad 4AC-Q^2\geq\theta.$$
\end{lem}

\begin{proof}
%{Proof of Lemma \ref{lem:sub}}
Let $A$, $Q$, $C$ be as in the statement. Set $Lu:= Au''-Qu'+Cu$. For $\tau, k >
0$, it
holds that
$$(1 -\tau)4A(C - k)-Q^2 = (1 - \tau)(4AC -Q^2 - 4Ak)- \tau Q^2 \geq (1 -
\tau)\theta - \tau\sigma^2 - 4\beta k.$$
Hence, it is possible to choose $\tau$, $k$, only depending on $\beta$, $\sigma$,
$\theta$, in such a way that
$4A(C - k) - Q^2 \geq 4\tau A(C - k)$. As a consequence, the function $g(\rho) :=
\rho^n$ satisfies
\[\begin{split}
(L - k)g &=
\left(A(n^2 - n) - Qn\rho + (C - k)\rho^2\right)
\rho^{n-2}\\
& \geq \left(\frac{ 4A(C - k) - Q^2}{4(C - k)}
n - A\right)n\rho^{n-2}\\
&\geq (\tau n - 1)An\rho^{n-2}.
\end{split}\]
There exists then $n \geq 3$ such that $(L-k)g \geq 0$ in $\R$. We define the
function $h$ and
the constant $r$ in the following way:
$$h(\rho) := \chi(\e\rho)g(\rho) + \e^{-n}(1 - \chi(\e\rho)), 
\qquad r := \e^{-1},$$
where $\e$ is a positive constant that will be chosen later and
$\chi$ is a smooth function satisfying
$$\chi = 1\text{ in }(-\infty, 1/2],\qquad
\chi'\leq0 \text{ in }(1/2, 1),\qquad 
\chi = 0 \text{ in } [1,+\infty).$$
By direct computation one sees that $h' > 0$
in $(0, r)$. For
$0 <\rho\leq (2\e)^{-1}$, it holds that $Lh=Lg > 0$. Let $(2\e)^{-1} <\rho<
\e^{-1}$.
Using the inequality $Lg \geq kg$, we get
\[\begin{split}
Lh &= \chi Lg + A\chi''\e^2g + 2A\chi'\e g' - Q\chi'\e g +
\e^{-n}[- A\chi''\e^2 + Q\chi'\e + C(1 - \chi)]\\
&\geq
\left[2^{-n}k\chi-2\beta|\chi''|\e^2+2\beta\chi'n\e^2+\sigma\chi'\e
+ C(1 - \chi)\right]\e^{-n}.
\end{split}\]
Notice that, by hypothesis, $C\geq\frac\theta{4\beta}$.
Hence, $2^{-n}k\chi+C(1 - \chi)\geq\min(2^{-n}k,\frac\theta{4\beta})$. As a
consequence,
$$Lh\geq\left[\min\left(2^{-n}k,\frac\theta{4\beta}\right)
-(2\beta|\chi''|\e-2\beta\chi'n\e-\sigma\chi')\e\right]\e^{-n}.$$
Therefore, for $\e>0$ small enough, $h$ satisfies
$Lh>0$ in $((2\e)^{-1}, \e^{-1})$, and then in $(0,r)$.
\end{proof}

%%%%%%%%%%%%%%%%%%%%%%%%%%%%%%%%%%%%%%%%%%%%
%%%%%%%%%%%%%%%%%%%%%%%%%%%%%%%%%%%%%%%%%%%%

\end{document}